\documentclass[12pt]{article}
\usepackage{amscd}
\usepackage{latexsym}
\begin{document}
\title{On open 3-manifolds proper homotopy equivalent to geometrically
simply-connected polyhedra  
\footnote{First version: {\tt June 8, 1998}.
 This version: {\today}.
This preprint is available electronically at 
          \tt  http://www-fourier.ujf-grenoble.fr/\~{ }funar}}
\author{
\begin{tabular}{cc}
 L. Funar &  T.L. Thickstun\\
\small \em Institut Fourier BP 74 
&\small \em Department of Mathematics \\
\small \em University of Grenoble I &\small \em Southwest Texas State
University\\
\small \em 38402 Saint-Martin-d'H\`eres cedex, France
&\small \em San Marcos, TX 78666, USA\\
\small \em e-mail: {\tt funar@fourier.ujf-grenoble.fr}
& \small \em e-mail: {\tt tt04@swt.edu } \\
\end{tabular}
}

\maketitle


\newcommand{\lam}{\lambda}

\newcommand{\vp}{\varphi}

\newcommand{\bd}{\partial}

\newcommand{\nul}{\emptyset}

\newcommand{\vs}{\vspace{\baselineskip}}

\newcommand{\Q}{{\bf Q}}

\newcommand{\Z}{{\bf Z}}


\newtheorem{df}{Definition}[section]

\newtheorem{lemma}{Lemma}[section]

\newtheorem{prop}{Proposition}[section]

\newtheorem{theo}{Theorem}[section]

\newtheorem{rem}{Remark}[section]

\newtheorem{cor}{Corollary}[section]


\begin{abstract}
We prove that an open 3-manifold proper homotopy equivalent to 
a geometrically simply connected polyhedron is simply connected at
infinity thereby  generalizing the
theorem proved by Po\'enaru in \cite{Po1}. 
\end{abstract}

\vspace{0.3cm}
{\em AMS MOS Subj.Classification}(1991): 57 M 50, 57 M 10, 57 M 30. 

\vspace{0.3cm}
{\em Keywords and phrases}: $H_3$-semi-dominated, simple-ended,
simply connected 3-manifold, geometric simple connectivity, 
simple connectivity at infinity.

\section{Introduction}
The immediate antecedent to 
this paper is \cite{Po1}, the principal 
theorem of which is the following.
\begin{theo}
(V.Po\'enaru)
If $U$ is an open simply connected  3-manifold and, for some $n$, 
$U\times D^n$ has a handlebody decomposition without 1-handles 
then $U$ is geometrically simply connected hence simply connected at infinity. 
\end{theo}

\noindent {\bf Note}: $D^n$ denotes the $n$-ball; see \cite{Po1} for the definition
of ``handlebody decomposition without 1-handles ``; all 3-manifolds we consider in the sequel will be orientable
unless the contrary is explicitly stated. 

A non-compact polyhedron  $U$ is {\em simply connected at infinity} 
(s.c.i.), and we write also  $\pi_1^{\infty}(U)=0$, if given   
a compactum (i.e. a compact set) $X\subset U$ there exists another compactum $Y$ with 
$X\subset Y\subset U^3$, such that any loop in $U-Y$ is null-homotopic in 
$U-X$. (Some authors call this $\pi_1$-{\it triviality at 
infinity} or {\it 1-LC at infinity} and reserve the term 
s.c.i. for the special case in which $Y$ can be chosen so that, in
addition, $U-Y$ is connected. These notions are equivalent for one
ended spaces, such as contractible spaces. 

A non-compact polyhedron  $P$ is  
 {\em geometrically simply connected (g.s.c.)}  if it can be exhausted by
 compact 1-connected subpolyhedra. 
It is easily seen that any non-compact polyhedral manifold $U$ which
has a handlebody decomposition with no 1-handles is  g.s.c.   
In addition the projection map $p:U\times D^n\to U$ is a proper
simple-homotopy equivalence (defined in \cite{Sie}). 
In \cite{Po1}
Po\'enaru hinted at the conjecture which results when the hypothesis 
of the above-stated theorem is replaced by the (therefore) weaker
hypothesis that $U$ be proper simple-homotopy equivalent to a g.s.c. 
polyhedron. This conjecture was subsequently established in \cite{Fun}
using the techniques of \cite{Po1}. 
The following theorem is an immediate corollary of the principal 
result of this paper (which is proven using only basic 3-manifold
theory). It further generalizes the theorem stated above. 

\begin{theo} 
Any open 3-manifold which is proper homotopy equivalent to a
geometrically simply connected  polyhedron is simply-connected at infinity. 
\end{theo}

\begin{rem}
If $M^n$ ($n>3$) is a compact, contractible $n$-manifold with 
non-simply-connected boundary (e.g. those constructed in 
\cite{Mazur0} and \cite{Po}) then $int(M^n)$ is easily seen to be 
g.s.c. but not s.c.i. This demonstrates that
the above theorem cannot be extended to include open contractible 
manifolds of dimension greater than three. 
Notice that for simply connected open 3-manifolds the g.s.c. and 
the s.c.i. are equivalent. 
\end{rem}

\noindent {\bf Provisos}: We remain in the polyhedral category
throughout and all homology groups are with ${\bf Z}$ coefficients. 

\vspace{0.5cm}

\noindent {\bf Acknowledgements}: We are indebted to the referee 
for his valuable suggestions which simplified some proofs and improved
the exposition. The first author wishes to
express his gratitude to V.Po\'enaru for the discussions they had on 
this subject.  

Several months after the appearance of \cite{Fun} in 
preprint form, the second-named author of this paper communicated 
to the first-named author an outline of the proof of the principal
result found here (the details were then worked out jointly).

\section{Statement of results}
We first require the following definitions. 
\begin{df}{\rm 
A proper map $f:X\to Y$ is {\em  $H_3$-nontrivial}  
if given  non-null compacta $L\subset Y$ and 
$K\subset X$ such that $f(X-K)\subset Y-L$ 
then $f_*:H_3(X,X-K)\to H_3(Y,Y-L)$ is nontrivial 
(i.e. its image is not a singleton).} 
\end{df}

\begin{df}{\rm 
Given noncompact polyhedra $X$ and $Y$ we say {\em $Y$ is 
$H_3$-semi-dominated by $X$} if there exists an $H_3$-nontrivial 
proper map $f:X\to Y$. }
\end{df}

\begin{df} {\rm 
An open connected 3-manifold $U^3$ is {\em simple-ended }if 
it has an exhaustion $\{M_i\}_{i=1}^{\infty}$ by 
compact 3-submanifolds where, for all $i$, the genus  
of $\partial M_i$ is zero. }
\end{df}

\begin{rem}
\begin{enumerate}
\item If $U$ is an open orientable 3-manifold and $Y\subset U$ 
a non-null compactum then $H_3(U,U-Y)$ is nontrivial. 
\item 
If $f,g:X\to Y$ are properly homotopic maps where $X$ and $Y$ are
open, connected 3-manifolds and $f$ is 
$H_3$-nontrivial then  $g$ is $H_3$-nontrivial. 
 
\end{enumerate}
\end{rem}
{\em Proof of 1.:} Since we are in
the polyhedral category $Y$ is a compact polyhedron and so has a
regular neighborhood $N$. Since $N-Y$ deformation retracts onto 
$\partial N$ we have that $H_3(U,U-Y)$ is isomorphic to 
$H_3(U,U-int\;N)$. Taking a triangulation of the pair $(U,N)$ and
using the orientation of $U$ we see that $N$ is a relative cycle
representing a nontrivial element of $H_3(U,U-int\;N)$. Moreover, this 
is a free abelian group with basis elements represented by the
components of $N$. This proves the first remark. 

{\em Proof of 2.:} Suppose that $K$, $L$ are compacta 
in $X, Y$, respectively, with $g(X-K)\subset Y-L$. Then $K\supset
g^{-1}(L)$. Let $F:X\times I\to Y$ be a proper homotopy with 
$F(x,0)=f(x)$ and $F(x,1)=g(x)$. Since $F^{-1}(L)$ is compact, there 
is a connected compactum $M$ of $X$ such that $M \times I \supset F^{-1}(L) \cup (K \times \{1\})$. Let $N$
be a regular
neighborhood of $M$ in $X$. Then $N$ represents a generator of the
infinite cyclic
group $H_3(X,X-int \, N)$ which is isomorphic to $H_3(X,X-M)$. Now
$f^{-1}(L) \subset
M$ implies that $f(X-M) \subset Y-L$. Since $f$ is $H_3$-nontrivial we
have that
$f(N)$ represents a nontrivial element of $H_3(Y,Y-L)$. Since $K \subset
int \, N$ we
have that $N$ represents an element of $H_3(X,X-K)$. Since $F^{-1}(L)
\subset int \, N
\times I$ we have that $g(N)$ is homologous to $f(N)$ in $(Y,Y-L)$, and
so $g$ is
$H_3$-nontrivial. $\Box$

Our main result is the following.
\begin{theo}
An open, connected, orientable 3-manifold $U$ is 
$H_3$-semi-dominated by a g.s.c. polyhedron if and only if $U$ is the 
connect-sum of a 1-connected, simply connected at infinity open
3-manifold and a closed, orientable 
3-manifold with finite fundamental group. 
\end{theo}

\begin{rem}
\begin{enumerate}
\item Theorem 1.2. is a
  corollary of Theorem 2.1. In fact, let $f:P
\rightarrow U$
be a proper homotopy equivalence with proper homotopy inverse $g:U
\rightarrow P$.
Suppose $K$, $L$ are compacta in $P$, $U$, respectively, such that
$f(P-K) \subset
U-L$. Let $M=g^{-1}(K)$. Then $g(U-M) \subset P-K$. Since $f \circ g$ is
proper
homotopic to $id_U$ and $id_U$ is $H_3$-nontrivial, we have by Remark
2.1 (2) that
$(f \circ g)_*$ is nontrivial and hence $f_*$ is nontrivial.
\item The proof of the ``if'' part of the theorem is very brief. Just
  observe that if $U$ satisfies the hypothesis then the universal 
covering of $U$ is 1-connected and simply connected at infinity 
(hence g.s.c. -see (3) below) and the covering projection is proper
and has non-zero degree (hence is $H_3$-nontrivial).
(In the sequel when we refer to the hypothesis or conclusion of the
theorem we will mean the ``only if'' part.)
\item By the methods of \cite{Wa} the class of 1-connected, simply
  connected at infinity open 3-manifolds 
is  equal to each of the following two classes of open
3-manifolds. 
Those which can be constructed as follows: delete a tame,
0-dimensional, compact subspace from $S^3$, denote the result by $U$
and replace each element of a pairwise disjoint, proper family of
3-balls in $U$ by a 
homotopy 3-ball. 
Those each of which has an exhaustion $\{M_i\}_{i=1}^{\infty}$ by
compact 3-submanifolds where, for each $i$, $M_i$ is 
1-connected (and hence the genus of $\partial M_i$ is $0$)
and each component of $M_{i+1}-M_i$ is homeomorphic to a space obtained
by taking finitely many pairwise disjoint 3-balls in $S^3$, replacing
one 
by a homotopy 3-ball and deleting the interiors of the rest. 
\end{enumerate}
\end{rem}

To establish the theorem we will demonstrate the following three
propositions. 

\begin{prop}
If $U$ is as in the hypothesis of the theorem then $U$ is
simple-ended. 
\end{prop}

\begin{prop}
If $U$ is as in the hypothesis of the theorem then $\pi_1(U)$ is a
torsion group. 
\end{prop}

\begin{prop}
If $U$ is an open, connected simple-ended 3-manifold such that 
$\pi_1(U)$ is a torsion group then $U$ is as in the conclusion of the
theorem. 
\end{prop}

\section{Proof of Proposition 2.1}
\begin{lemma}
Suppose the following: $U$ is an open, orientable 3-manifold; 
$K$ is a compact, connected 3-submanifold of $U$ such that each
component of $cl(U-K)$ is noncompact and has connected boundary; and 
$f:(M,\partial M)\to (U,U-K)$ is a map of a compact connected
3-manifold with boundary such
that $\partial M$ has genus zero and 
$f_*:H_3(M,\partial M)\to H_3(U,U-K)$ is nontrivial. 
Then there exists a compact, connected 3-submanifold $N$ of $U$ such that
$K$ is in $N$ and the genus of $\partial N$ is zero. 
\end{lemma}
{\it Proof:} 
It will suffice to find, for each component $V$ of $cl(U-K)$ a 
compact, connected 3-submanifold $N(V)$ of $V$ such 
that $\partial V$ is in $N(V)$ and the genus of $\partial
N(V)-\partial V$ is zero. So let $V$ be such a component. We assume
$f$ is transverse to $\partial V$ and denote the intersection of 
$f^{-1}(V)$ and $\partial M$ by $C$. The
orientation of
$U$ determines an orientation of $V$, which in turn determines an
orientation of $\bd V$.
Similarly the orientation of $M$ determines an orientation of the
3-dimensional
submanifold $f^{-1}(V)$, which in turn determines an orientation of
$f^{-1}(\bd V)
\subset \bd f^{-1}(V)$.

From the hypothesis on $f_*$ we derive 
that
$(f|_{f^{-1}(\bd V)})_*$ carries the homology class of $f^{-1}(\bd V)$,
as
determined by the orientation above, to a nontrivial element of $H_2(\bd
V)$. 
In fact, consider the commutative diagram
\[ \begin{CD} H_3(M,\bd M) @>{\bd_M}>> H_2(f^{-1}(\bd V)) \\
              @V{f_*}VV              @VV{(f|_{f^{-1}(\bd V})_*}V \\
              H_3(U,U-int\,K) @>{\bd_U}>> H_2(\bd V) \end{CD} \]
where $\bd_M$ comes from the relative Mayer-Vietoris sequence of 
$(M,\bd M)=(f^{-1}(V), C) \cup (\overline{M-f^{-1}(V)}, \bd M-C)$
and $\bd_U$ comes from the relative Mayer-Vietoris sequence of
$(U,U-int \, K)=(V,V) \cup (\overline{U-V},\overline{U-(V \cup K)})$.
A chase through the zig-zag lemma shows that $\bd_M$ takes the
orientation class of $(M, \bd M)$ to the class of $f^{-1}(\bd V)$
and $\bd_U$ takes the generator of the infinite cyclic group
$H_3(U,U-int \, K)$ to the class of $\bd V$. The proof is completed
by noting that inclusion induces an isomorphism of $H_3(U,U-int \, K)$
with $H_3(U,U-K)$ (since $U-int \, K$ and $U-K$ both deformation
retract onto $U$ minus the interior of a regular neighborhood of $K$).

Note that 
$f{\mid}_C:C\to V$ and  $f\mid_{f^{-1}(\partial V)}:f^{-1}(\partial
V)\to V$ are homologous and hence $f{\mid}_C$ is homologous (in $V$)
to a nonzero multiple of $\partial V$. 

Applying the prime factorization theorem for compact 3-manifolds to a 
regular neighborhood of $f(C)$ in $int(V)$ we conclude the existence 
of an embedding $R$ in $int(V)$ where $R$ is a closed, oriented 
surface of genus zero such that $R$ is homologous to $f{\mid}_C$ in
$int(V)$
(recall that the prime factorization  and the sphere theorem imply
that $\pi_2$ 
of a compact 3-manifold is generated, as a $\pi_1$-module, by a finite
family of pairwise disjoint embedded 2-spheres, and by the Hurewicz
theorem $H_2$ is isomorphic to $\pi_2$ modulo the action of $\pi_1$). Now let $W$ be a
regular neighborhood of $R$ in $int(V)$ and denote by $N(V)$ that
component of $cl(V-W)$ containing $\partial V$. It remains only to
show that $N(V)$ is compact. Suppose otherwise. Then there 
exists a proper ray in $V$ extending from $\partial V$ to infinity 
and avoiding $R$. Such a ray has intersection number one with
$\partial V$ but intersection number zero with $R$. This contradicts
the fact that $R$ is homologous in $V$ to a multiple of $\partial V$. 
 $\Box$

\begin{df}\label{adm}
{\rm An {\em admissible pair} is a map 
$f:(X,Y)\to (M,M-L)$ and a subspace $K\subset L$ satisfying the
following conditions:
\begin{enumerate}
\item $(X,Y)$ is a pair of compact simplicial complexes and 
$X$ is simply connected. 
\item $M$ is a compact orientable 3-manifold and $K$ and $L$ are 
compact 3-submanifolds with $K\subset int(L)$, $L\subset int(M)$ and 
$L$ connected. 
\item The map $f$ is simplicial and $M$ is an (abstract) regular
  neighborhood of $f(X)$. 
\item Only one component of $X-f^{-1}(\partial L)$ has image under $f$ 
meeting $K$. 
\item The map $f_*:H_3(X,Y)\to H_3(M,M-L)$ is non-trivial. 
\end{enumerate}
We will refer to the $M$ above as the {\em target} of the admissible pair. }
\end{df}
\noindent {\bf Notation}: If $X\subset P$  then $\partial X$ denotes 
the frontier of $X$ in $P$. 
\begin{lemma}\label{dehn}
If $f:P\to U$ is the map of the theorem and $K$ is a compact
3-submanifold of $U$ then we can choose $X\subset P$ and 
a compact 3-submanifold $L$ of $U$ such that the pair 
$\{f|_X:(X,\partial X)\to (M,M-L), K\}$ is admissible where 
$M$ is a regular neighborhood of $f(X)$ in $U$. 
\end{lemma}
{\it Proof:} 
Consider  $L_0$ be a regular neighborhood of $K$ in $U$. 
The hypothesis implies the existence of 
some $X_0$ with
 $f^{-1}(L_0)\subset X_0\subset P$, such that 
$f_*:H_3(P,P-X_0)\to H_3(U,U-L_0)$ is non-trivial. 
Moreover once such an $X_0$ is chosen then larger 
$X$ with  $int(X)\supset X_0$ are also convenient for this purpose.
Moreover $f$ is proper implies the existence of an 
$X_1\supset X_0$  such that $f(\partial X_1)\cap
L_0=\emptyset$. Using excision on both sides above we find that 
the map 
$(f|_{X_1})_*:H_3(X_1,\partial X_1)\to H_3(M,M-L_0)$ is non-trivial,
where $M$ denotes a regular neighborhood of $f(X_1)$ in $U$. 
Denote the set $X_1$  with these properties by $X(L_0)$. 
Note that  $X(L_0)$ is defined  for any compact $L_0$ engulfing $K$. 
If $X-f^{-1}(\partial L_0)$ has at least two components 
each of which has  image meeting $K$ then  $f^{-1}(K)\cap
(X(L_0)-f^{-1}(\partial L_0))$ is not connected. 
Since $f$ is simplicial the latter is a simplicial subcomplex of $P$
and so it cannot be path connected.
In particular there exist points $x,y\in f^{-1}(K)$ such that 
any path connecting them in $X(L_0)$ should meet $f^{-1}(\partial
L_0)$. But $f^{-1}(K)$ is a compact (thus finite) simplicial 
complex because $f$ is proper hence it has  a finite number of
components. Consider some arcs joining these components in $P$. 
The union of these arcs with the components of $f^{-1}(K)$ is contained in some compact subset $K'\subset P$.
Consider now $L$ large enough such that 
$f^{-1}(\partial L)\cap K'= \emptyset$, and 
$K'\subset int(f^{-1}(L))$. Then $X(L)\supset f^{-1}(L)\supset K'$, and 
we claim that this $X$ fulfills all conditions needed. 
If $x,y$ are two points in $f^{-1}(K)$ then there exists a path
connecting them inside $K'$ and so there exists a path inside
$X(L)-f^{-1}(\partial L)$. $\Box$ 

\begin{lemma}
Given an open, connected 3-manifold $U$ and compactum $X$ in $U$ there
exists a compact 3-submanifold $K$ of $U$ containing $X$  such that
each component of $cl(U-K)$ is noncompact and has connected boundary. 
\end{lemma}
{\it Proof}: Let $M$ be a compact 3-submanifold of $U$ containing
$X$. 
Let $N$ be the union of $M$ and all compact components of
$cl(U-M)$. To obtain $K$ from $N$ add 1-handles to $N$ (in $cl(U-N)$)
which connect different components of $\partial N$ which are in the
same component of $cl(U-N)$. $\Box$ 

The proof of Proposition 2.1 will proceed by applying the tower
construction to the admissible pair of Lemma 3.2 (where $K$ is also 
chosen to satisfy the conclusion of Lemma 3.3) to obtain a map 
satisfying the hypothesis of Lemma 3.1. It is convenient to state
first the following definition.

\begin{df}\label{re}
{\rm A {\em reduction of the admissible pair} $\{ f_0:(X,Y)\to(M_0,M_0-L_0),
K_0\}$ is a second admissible pair 
$\{ f_1:(X,Y)\to(M_1,M_1-L_1), K_1\}$
such that there exists a map $p:M_1\to M_0$, the ``projection map'' 
of the reduction,  satisfying the following
conditions:
\begin{enumerate}
\item $p\circ f_1 = f_0$.
\item $p(K_1)=K_0$, $p(L_1)=L_0$, and the maps 
$p|_{L_1}:L_1\to L_0$, $p|_{K_1}:K_1\to K_0$ are boundary preserving. 
\item $p|_{L_1}:L_1\to L_0$ has non-zero degree. 
\item The complexity of $f_1$ is strictly less than the complexity of
  $f_0$ (where the complexity of a simplicial map $g$ with compact domain 
is the number of simplexes $s$ in the domain 
for which $g^{-1}(g(s))\neq s$). 
\item The image under $p$ of only one component of $M_1-\partial L_1$ meets
  $K_0$. Observe that, by (2), the image of $L_1$ must meet $K_0$. 
\end{enumerate}
}
\end{df}
\begin{lemma}\label{red}
If $f_1$ is a reduction of $f_0$ and $f_2$ is a reduction of $f_1$
then $f_2$ is a reduction of $f_0$. 
\end{lemma}
{\it Proof:} This is obvious.$\Box$ 
\begin{lemma}
Any admissible pair with non-simply connected target has a
reduction with simply-connected target. 
\end{lemma}
{\it Proof:} It will suffice to show that any admissible pair with
non-simply connected target has a reduction. Because then, by
iteration (which could occur at most finitely many times by condition
4) and applying Lemma 3.3 we obtain a reduction with simply connected
target. 

Let $\{f|_X:(X,Y)\to (M,M-L), K\}$ be the admissible pair.
Let $\widetilde{M}$ be the universal covering of $M$ and $p$ be the
covering projection. 
Since  $X$ is simply connected there exists 
a lift  $\widetilde{f}:X\to \widetilde{M}$ of $f$. 
We will show that 
$\{f_1:(X,Y)\to (M_1,M_1-L_1), K_1\}$
is a reduction, where 
\begin{enumerate}
\item $f_1=\widetilde{f}$.
\item $M_1$ is a regular neighborhood of $\widetilde{f}(X)$. 
\item $L_1$ is the only component of $p^{-1}(L)$ 
which is contained in $f_1(X)=\widetilde{f}(X)$. 
\item $K_1=p^{-1}(K_0)\cap L_1$. 
\item The covering map restricted to $M_1$ is the map $p$ from the
  definition \ref{re}. 
\end{enumerate}
Let us show that $f_1$ is well-defined and admissible. 

\vspace{0.5cm}
{\it $f_1$ is well-defined:}
First we show that $f_1(Y)\subset M_1-L_1$. We know that 
$f(Y)\subset M-L$ since $f$ is admissible, and so 
$f_1(Y)=M_1\cap p^{-1}(f(Y)) \subset M_1\cap p^{-1}(M-L)= 
M_1 -p^{-1}(L) \subset M_1-L_1$.

In order to have a consistent definition 
of $L_1$ we must show that there exists one and only one 
component of $p^{-1}(L)$ contained in $f_1(X)$. 

First we prove the existence. 
Since $f$ is non-trivial on $H_3$ (the condition (5) for $f$) we derive that 
$\widetilde{f}_*:H_3(X,Y)\to H_3(\widetilde{M}, \widetilde{M}-p^{-1}(L))$ 
is non-trivial. In particular the abelian group 
$H_3(\widetilde{M}, \widetilde{M}-p^{-1}(L))$ is non-zero. Since this 
group is  freely generated 
by an equivariant regular neighborhood of the compact components 
of $p^{-1}(L)$,
there  exists at least one such (
the deck transformations act transitively on the components 
of $p^{-1}(L)$ and so every component is compact). 
Since $L$ is connected $H_3(M,M-L)$ is generated 
by the orientation class of a regular neighborhood of $L$ and $f_*$ nontrivial implies that 
$L$ is in $f(X)$. 

Observe now that $p(M_1\cap p^{-1}(L))=L$ since $p(M_1)=M$, and also 
$p(M_1-M_1\cap p^{-1}(L))=M-L$. 
Then since $f$ is $H_3$ nontrivial the map 
$\widetilde{f}_*:H_3(X,Y)\to H_3(M_1,M_1-M_1\cap p^{-1}(L))$ should be 
non-trivial. The same argument used above shows that 
$\widetilde{f}(X)\supset M_1\cap p^{-1}(L)$. 
Let $L_1$ be a component of $p^{-1}(L)$ which meets $\widetilde{f}(X)$. 
Suppose that $L_1$ is not completely contained in $\widetilde{f}(X)$. 
Since $L_1$ is connected then $L_1 \cap (M_1 - \widetilde{f}(X))\neq
\emptyset$, or in other words the regular neighborhood of
$\widetilde{f}(X)$ meets a larger subset of $L_1$ than the image $\widetilde{f}(X)$. 
This contradicts the fact that $M_1\cap L_1 \subset \widetilde{f}(X)$. 
Thus any component of $p^{-1}(L)$ meeting $\widetilde{f}(X)$ is entirely
contained in $\widetilde{f}(X)$. At least one component has non-void
intersection with the image because 
$p(\widetilde{f}(X) \cap p^{-1}(L))=p(M_1\cap p^{-1}(L))= L$. 

Suppose now that there are two components, $L_1$ and $L_1'$ 
meeting $\widetilde{f}(X)$. The two components are then 
disjoint and contained in $\widetilde{f}(X)$. Furthermore there are two 
components of $M_1-p^{-1}(\partial L)$, namely $int(L_1)$ and 
$int(L_1')$ whose images under $p$ meet $K$. 
However there is only one component, say $\xi$, of $X-f^{-1}(\partial L)$ 
whose image by $f$ meets $K$. Then   
$\widetilde{f}(\xi)\subset M_1-\partial L_1 \cup \partial L_1'$ 
since $\widetilde{f}(\xi)$ 
avoids $p^{-1}(\partial L)$ and 
$\widetilde{f}(\xi)\supset int(L_1)\cup int(L_1')\supset K_1 \cup K_1'$
because $p(f(\xi))$ meets $K$. This is a contradiction as $\widetilde{f}(\xi)$ must
be connected.

\vspace{0.5cm}
{\it $f_1$ is admissible:}
Conditions (1-3) from definition \ref{adm} are immediate. 
The condition (5) for $f_1$ is 
satisfied since $f=p\circ f_1$ is $H_3$-nontrivial.
Finally (4) is implied by $f_1^{-1}(\partial L_1)=f^{-1}(\partial L)$.

\vspace{0.5cm}
{\it  $f_1$ is a reduction of $f$:}
With the exception of (3) and (5) we leave them to the reader. 
Condition (3) follows from the fact that $p|_{L_1}:L_1\to L$ 
is a covering map, and it is well-known that any covering map from one
compact, 
connected, orientable 3-manifold to another has non-zero degree. 
An easy argument is the following. Triangulate and orient the base space.
Lift the
triangulation and the orientation of each 3-simplex to the covering
space. If this
does not give an orientation of the covering space, then there will be a
2-simplex
which is a face of two 3-simplexes which are mapped to the same
3-simplex in the
base space, contradicting the fact that the map is a covering map. It
then follows
that a fundamental cycle for the covering space is sent to $n$ times the
fundamental
cycle for the base space, where $n$ is the number of sheets of the
covering.

We already saw before that
$int(L_1)$ is the only one component of $M_1-p^{-1}(\partial L)$ 
whose image meets $K$ hence establishing (5). 
$\Box$

{\it Proof of Proposition 2.1  from the Lemmas:}
Begin with the admissible pair of Lemma 3.2 where $K$ also satisfies
the conclusion of Lemma 3.3. Apply Lemma 3.5 to that admissible pair.
Note that the projection map of that reduction satisfies the
hypothesis of Lemma 3.1 whose application then completes the
proof. $\Box$

\section{Proof of Proposition 2.2}

\begin{lemma}
Suppose the diagram 
\[
\begin{array}{rcl}
M         & \stackrel{\varphi}{\longrightarrow} &  X    \\
h\downarrow&  &  \downarrow \psi \\
N       & \subset & U 
\end{array}
\]
is commutative and satisfies the following conditions:
\begin{enumerate}
\item
$M$ and $N$ are compact, connected,  orientable 3-manifolds with 
$h(\partial M)\subset \partial N$. 
\item $h$ has non-zero degree. 
\item $X$ is simply-connected. 
\end{enumerate}
Then $e_{\sharp}(\pi_1(N))$ is a torsion subgroup of $\pi_1(U)$ 
(where $e$ denotes the inclusion $N\subset U$).
\end{lemma}
{\em Proof:} Let
$p:\widetilde{N} \rightarrow
N$ be the covering map such that
$p_{\sharp}(\pi_1(\widetilde{N}))=h_{\sharp}(\pi_1(M))$.
Then $h$ lifts to $\tilde{h}:M \rightarrow \widetilde{N}$. Since the
degree of $h$ is
non-zero $H_3(\widetilde{N},\bd\widetilde{N})\neq 0$, and so
$\widetilde{N}$ is compact
and thus $p$ is finite sheeted. Hence $h_{\sharp}(\pi_1(M))$ has finite
index in
$\pi_1(N)$. (This is a standard argument.) Thus for any element of
$\pi_1(N)$ some
power will be in the image of $h_{\sharp}$ and therefore will be trivial
in $\pi_1(U)$.
$\Box$

Now to prove Proposition 2.2 it will suffice to show that if $N$ is a
compact, connected 3-submanifold of $U$ then $e_{\sharp}(\pi_1(N))$ 
is a torsion subgroup of $\pi_1(U)$. Let $f:X\to U$ 
be as in the definition of ``$H_3$-semi-dominated'' (where $X$ is
g.s.c.).
By excision we have 
$f_*:H_3(Y,\partial Y)\to H_3(N,\partial N)$ is nontrivial, where 
$Y=f^{-1}(N)$. We can ``realize'' 
any element of $H_3(Y,\partial Y)$ by a map 
$g:(M,\partial M)\to (Y,\partial Y)$, where $M$ is a compact 
orientable 3-manifold (i.e. $g_*:H_3(M,\partial M)\to H_3(Y,\partial
Y)$ sends the orientation class of $M$  to the preassigned element of 
$H_3(Y,\partial Y)$ - see \cite{Thom}). The 3-manifold $M$
coming from
Thom's theorem might not be connected. But since $(f \circ g)_*:
H_3(M,\bd M)
\rightarrow H_3(N,\bd N)$ is non-trivial there will be some component
$M_0$ of
$M$ such that $((f \circ g)|_{M_0})_*:H_3(M_0,\bd M_0) \rightarrow
H_3(N, \bd N)$
is nontrivial. Now apply the lemma with $h=f \circ (g|_{M_0})$. $\Box$

\section{Proof of Proposition 2.3}
\begin{lemma}
If $A$ and $B$ are groups and $a\in A$, $b\in B$ are neither the
identity then $a*b$ is not a torsion element in $A*B$ (the free
product of $A$ and $B$). 
\end{lemma}
{\em Proof:} This is a standard fact from combinatorial
group theory.
Every element of finite order in $A*B$ is conjugate to an element of $A$
or of $B$
(Corollary 4.1.4 of \cite{MKS}), and every element of $A*B$ is conjugate
to a
cyclically reduced word which is unique up to cyclic permutation
(Theorem 4.2 of
\cite{MKS}). Since $a*b$ is cyclically reduced it cannot be conjugate to
an element
of $A$ or of $B$ and hence cannot have finite order. $\Box$ 

\begin{cor}
The prime factorization of a closed 3-manifold whose fundamental group
is a torsion group can have at most one non-simply connected factor. 
\end{cor}

Now to prove Proposition 2.3 we assume that $U$ is not simply-connected 
(otherwise we are done) and let $M$ be a non-simply connected
3-submanifold of $U$ such that $\partial M$ has genus zero. 
We can express $M$ as a connect-sum of a punctured homotopy 3-ball 
and a closed orientable 3-manifold $N$ where $\pi_1(N)$ is a torsion
group. 
By the above corollary we can assume $N$ is irreducible. 
By \cite{Eps} any orientable, irreducible closed 
3-manifold with torsion must have finite fundamental group. It remains
only to show that $U$ can have no other non-simply connected factor
but this also follows from the Corollary. $\Box$

\begin{small}
\bibliographystyle{plain}

\end{small}

\end{document}